\documentclass[11pt]{article}
\usepackage{amscd,amssymb,amsmath,verbatim,amsfonts,amsthm,color}
\newtheorem{theorem}[equation]{Theorem}

\newtheorem{lemma}[equation]{Lemma}

\numberwithin{equation}{section}
\newtheorem{defn}[equation]{Definition}
\date{}

\newcommand{\im}{\operatorname{Im}}

\newcommand{\reals}{{\mathbb R}}

\newcommand{\sphere}{{\bf S}}

\newcommand{\RR}{\mathbb{R}}
\newcommand{\HH}{\mathbb{H}}

\newcommand{\calC}{\mathcal C}
\newcommand{\calP}{\mathcal P}
\newcommand{\calS}{\mathcal S}
\newcommand{\calE}{\mathcal E}
\newcommand{\calD}{\mathcal D}

\newcommand{\del}{\partial}
\newcommand{\e}{\epsilon}

\usepackage{enumerate}
\newcommand{\pd}[1][]{\partial_{#1}}
\newcommand{\PD}[1][]{D_{#1}}
\newcommand{\norm}[2][]{\left\| #2 \right\|_{#1}}

\newcommand{\grad}{\nabla}
\newcommand{\lap}{\Delta}
\DeclareMathOperator{\LP}{LP}
\DeclareMathOperator{\supp}{supp}

\newcommand{\CC}{\mathbb C}

\newcommand{\dt}{\;dt}

\newcommand{\dz}{\;dz}
\newcommand{\dsigma}{\;d\sigma}
\newcommand{\dtheta}{\;d\theta}

\newcommand{\coeff}{\varphi}
\newcommand{\wf}{\alpha}
\newcommand{\ph}{\Gamma}
\newcommand{\phr}{S} 

\newcommand{\calH}{\mathcal H}

\begin{document}

\title{Some global aspects of linear wave equations} 
\author{Dean Baskin \\ Northwestern University \and Rafe Mazzeo \\ Stanford University}

\maketitle

\begin{abstract}
This paper surveys a few aspects of the global theory of wave equations. This material is structured
around the contents of a minicourse given by the second author during the CMI/ETH Summer
School on evolution equations during the Summer of 2008.
\end{abstract}

\section{Introduction}
The week-long minicourse on which this brief survey paper is based came after a vigorous, detailed and outstanding 
series of lectures by Jared Wunsch on the applications of microlocal analysis to the study of linear wave equations. 
Both lecture series took place at the Clay Mathematics Institute Summer School at ETH Z\"urich in 2008. The goal of 
this minicourse was to describe a few topics which involve global aspects of wave theory, relying at least to some 
extent on the microlocal underpinnings from Wunsch's lectures.  The first of these topics is an account of 
some striking consequences that can be derived from the finite propagation speed property. While this had
been applied in various interesting ways before, the systematic development of this principle appears in 
the very influential paper of Cheeger, Gromov and Taylor \cite{CGT}.  We recall how this property, applied to
solutions of the wave equation associated to a Laplace-type operator, can be used to obtain estimates 
for solutions of various related operators. We present only one application of this, which is a lovely argument due to 
Gilles Carron which estimates the off-diagonal decay profile of the Green function for generalized Laplace-type
operators on globally symmetric spaces of noncompact type. This result had caught the lecturer's eye in the months 
before this Clay meeting and nicely illustrates the unexpected power of the finite propagation speed method.  
Following this, the remainder of the lectures reviewed several diifferent approaches to scattering theory and described
a few of the relationships between these.  The primary goal, however, was to introduce the Friedlander radiation 
fields and explain how they give a concrete realization of the Lax--Phillips translation representation. We follow
suit here, recalling the outlines of a few of the numerous successful approaches to scattering theory 
and culminating in a discussion of these radiation fields. 

This paper attempts to give some feel for what was presented in these lectures. The reader should be warned
that the topics covered here are in many places old-fashioned and we omit any mention of many of the 
most important recent advances and trends in scattering theory. The material here is meant to indicate a few
things that can be accomplished, often with not very sophisticated machinery by modern standards. We typically
make very restrictive assumptions in order to convey the main essence of the ideas. We give references
for further reading interspersed inter alia, but do not make any claim to a comprehensive bibliography.

The material assembled here is based on the notes of the first-named author; the lecturer (and second author) is
extremely grateful to him not only for this careful recording of the lectures, but also for his enthusiasm during
the lectures and his very substantial assistance in writing this paper.  We did discuss at some point, but later abandoned,
the possibility of writing a much more exhaustive treatment of some of the topics here, particularly the theory
of radiation fields. That will unfortunately have to wait for another day and other authors. We hope that this survey 
accomplishes what the original lectures also attempted, which is to whet the reader's curiosity to learn more 
about this subject. Needless to say, wave theory is an immense subject and we mention here only a very
small set of possible topics.

Throughout this paper we focus on properties of solutions, and of the solution operator, for the wave operator 
\begin{equation}
\Box_V = D_t^2 - L, \quad \mbox{where}\qquad L = \nabla^* \nabla + V
\label{defbox}
\end{equation}
acting on sections of some bundle $E$ over a Riemannian manifold $(M,g)$, where $\nabla$ is the covariant
derivative of some connection on $E$ and $V$ is a (self-adjoint) potential of order $0$, which can either be scalar 
or an endomorphism of $E$. For simplicity we typically assume that $V$ is smooth and compactly supported, although 
neither of these properties are present in almost any of the interesting physical or geometric applications. 
Furthermore, we often discuss only the scalar Laplacian and its perturbations, although the extension of
all results below to this slightly more general framework is usually just notational. Finally, here and below 
we write $\PD = \frac{1}{i} \pd$. 

As noted above, we take advantage of the luxury of being able to refer back to the excellent lecture notes
by Jared Wunsch \cite{W} covering his longer minicourse. Those notes provide a nice introduction for many 
central themes and results in the subject, including the existence of solutions of the equation $\Box u = f$ 
with vanishing Cauchy data, or of $\Box u = 0$ with prescribed nonzero Cauchy data, along a noncharacteristic 
hypersurface, the positive commutator method leading to H\"ormander's renowned theorem on propagation of 
singularities of solutions, the finite propagation speed property, and much else besides. Using this as a blanket 
resource, we can dive right into the material at hand. 

There are now many terrific monographs concerning the local and global aspects of wave equations. Michael Taylor's 
three-volume series \cite{Taylor} belongs high on this list; it contains an amazing amount of information about many 
different topics.  Other recent monographs with a particular focus on hyperbolic equations include those 
by Alinhac \cite{Al}, Lax \cite{Lax}, Rauch \cite{Rauch}; we mention also the new book by Zworski on
semiclassical analysis \cite{Zwo}. 

The authors are very grateful to the Clay Foundation for making this Summer School possible
-- it was a lot of fun and the large attendance and enthusiasm of the participants was amazing.
We also appreciate the forebearance by the editors of this volume for their (relative) tolerance for
the length of time between the original lectures and when this paper was finally written. 
Both authors are very grateful to many people for teaching us about many of the topics here.
We thank, in particular, Gilles Carron, Richard Melrose, Gunther Uhlmann, Andras Vasy and Jared Wunsch. 
Gilles Carron and Andras Vasy also gave some helpful remarks on this paper. 

D.B. is supported by NSF Postdoctoral Fellowship DMS-1103436; R.M. is supported by the NSF grant DMS-1105050.

\section{Finite propagation speed and its consequences}
Although \cite{W} contains a proof of the basic finite propagation speed property for the operator $\Box_V$, we begin by
recalling this familiar argument very briefly. We then show how using the functional calculus one can write the Schwartz kernels 
of various functions of the elliptic operator $L$ in terms of the Schwartz kernel of the wave operator. This leads directly
to the important Cheeger-Gromov-Taylor theory which uses finite propagation speed to obtain interesting estimates for 
these Schwartz kernels. We illustrate this with an outline of Carron's estimates for the resolvent and heat kernel of 
generalized Laplacians on symmetric spaces of noncompact type. 

\subsection{Finite propagation speed}
The fundamental identity behind finite propagation speed is the observation that for any sufficiently regular function $u$, 
\begin{equation}
\mathrm{div}_{x} (u_t \nabla u) = u_t \Box_0 u + \frac12 \del_t ( u_t^2 + |\nabla u|^2).
\label{1}
\end{equation} 
We suppose that the space on which we are doing calculations has a global time function $t$ and moreover, splits as 
$\RR \times M$, with a static Lorentzian metric $-dt^2 + h$, where $(M,h)$ is a Riemannian manifold. 
A hypersurface $Y \subset \RR \times M$ is called spacelike if its unit normal $\nu$ (with respect to this Lorentzian metric) 
satisfies $\nu \cdot \nu < 0$.  Suppose that $\Omega \subset \RR \times M$ is a domain bounded by two
spacelike hypersurfaces, $\del \Omega = Y_1 \cup Y_2$, which meet transversely along a codimension two submanifold,
and that $u$ is a solution of the homogeneous wave equation, $\Box_0 u = 0$. Integrate \eqref{1} over $\Omega$. The left 
side is transformed using the divergence theorem; the first term on the right vanishes while the second term is also transformed 
to a boundary integral. If $\nu_j = (\nu_{t,j}, \nu_{x,j})$ is the upward-pointing unit normal to $Y_j$, decomposed into its 
vertical ($t$) and horizontal ($x$) components, then we obtain
\begin{multline*}
\int_{Y_1} ( |u_t|^2 + |\nabla u|^2 ) |\nu_{t,1}| - 2 u_t \cdot \del_{\nu_1} u |\nu_{x,1}|  \\
=  \int_{Y_2} ( |u_t|^2 + |\nabla u|^2 ) |\nu_{t,2}| - 2 u_t \cdot \del_{\nu_2} u |\nu_{x,2}|.
\end{multline*}
Since $\nu_j$ is timelike, the integrand on each side is bounded from below by $c ( |u_t|^2 + |\nabla u|^2)$ for some $c > 0$
which depends on $Y_j$. We conclude that if $u_t = \nabla u = 0$ on $Y_1$, then these same quantities must also vanish 
on $Y_2$.  Finally, if $\Omega$ is foliated by spacelike hypersurfaces, then the vanishing of $(u_t, \nabla u)$ on the bottom 
(spacelike) boundary of $\Omega$ can be propagated throughout this entire region, and hence if $u$  vanishes at $Y_1$, 
then $u \equiv 0$ in $\Omega$. 

If we consider wave operators with terms of order $0$ or $1$, then this calculation can be adapted to show that if 
$\Omega$ is foliated by spacelike hypersurfaces $Z_s$, $0 \leq s \leq 1$, then the integral over $Z_s$ of 
$|u_t|^2 + |\nabla u|^2$ satisfies a differential inequality, and the fact that it vanishes when $s = 0$ implies that it 
vanishes for all $s \leq 1$. 

To interpret this calculation, we observe that there many natural domains $\Omega$ which can be foliated by spacelike 
hypersurfaces in this way. Indeed, suppose that $p = (t_1,x_1)$ is any point, and $\calD^-_{t_1,x_1}$ denotes the (backward) 
domain of dependence of this point, i.e. the set of points in $\RR \times M$ which can be reached by timelike paths 
travelling backward in $t$ and emanating from $(t_1,x_1)$. Let $Y$ be one of the level sets $\{t=t_0\}$ where $t_0 < t_1$.
Then the region $\Omega = \{(t,x) \in \calD^-_{t_1,x_1}: t \geq t_0\}$ can be shown to have a spacelike foliation
by submanifolds $Y_s$ which all intersect along the submanifold $\{(t,x) \in \calD^-_{t_1,x_1}: t = t_0\}$.  Thus any
homogeneous solution of $\Box_V u = 0$ which vanishes along with its normal derivative along $\{t = t_0\}$ vanishes
throughout this $\Omega$. This implies that if the Cauchy data of $u$ at $t_0$ is supported in some
subset $K$, then the Cauchy data of $u$ at $t_1 = t_0 + \tau$, where $\tau > 0$, is supported in the subset 
$K_\tau = \{(t_1,x): \mbox{dist}_g (x,K) \leq \tau\}$, which is precisely what is meant by saying that the support of a solution
propagates with speed $1$. For more general variable coefficient hyperbolic equations, the speed of propagation
may be variable but is still finite. 

\subsection{Cheeger-Gromov-Taylor theory}
Consider the fundamental solution for the problem 
\[
\Box_V u = 0, \quad \left. u \right|_{t=0} = \phi,\ \left. \del_t u \right|_{t=0} = 0.
\]
It is customary to write this solution operator as $\cos (t \sqrt{L})$, so that the solution $u(t,x)$ is equal to
$\cos (t\sqrt{L}) \phi$. We assume for simplicity that $L$ has no negative eigenvalues so there is no
ambiguity in this expression. This is an instance of the functional calculus for self-adjoint operators, which are 
defined in purely abstract terms using the spectral theorem, on one hand, but which can be used to 
describe solution operators for various equations involving $L$.  There are many interesting examples, including
prominently the resolvent and heat operator 
\[
R_L(\lambda) := (L - \lambda^2)^{-1} \qquad \mbox{and} \qquad e^{-tL}.
\]

The abstract definitions of these operators (i.e. defined using the spectral theorem) are all well and good, but
in order to use them one usually wishes to know much more about their mapping properties.  For example, 
a priori, using only these abstract definitions, we only know how one of these functions of $L$ acts on $L^2$
functions, but not on other function spaces. The goal then is to obtain a more concrete understanding of
the Schwartz kernels of any one of these operators.  Of course, there is a lot of theory devoted to doing just
this. Thus the classical theory of pseudodifferential operators gives a nice picture of the resolvent for any
one value of $\lambda$, while the theory of semiclassical pseudodifferential operators provides a means to
understand this family of operators as $\lambda$ tends to infinity in various directions in the complex plane. 
Similarly, the well-known heat-kernel parametrix construction,
cf.\ \cite{BGV}, gives a way to understand the
asymptotic behaviour of the Schwartz kernel of the solution operator for the heat equation in various regimes
of the space $\RR^+ \times M \times M$.  These theories and constructions give very precise information,
but are often very intricate, and furthermore, it is often hard to use these ideas directly to say anything 
interesting about global behaviour of these Schwartz kernels.   The idea in \cite{CGT} is that one can
extract, often in a rather simple way, some very useful global behaviour of these kernels using mainly
the finite propagation speed property of $\cos (t \sqrt{L})$ and some other simple properties, such
as the fact that the norm of $\cos (t\sqrt{L})$ as a bounded operator on $L^2$ never exceeds $1$. 

To explain this, suppose that $f(s)$ is a smooth, even function on $\RR$ which decays sufficiently rapidly so 
that the following manipulations are justified. Assuming $L \geq 0$ for simplicity, we define $f(\sqrt{L})$ using
the spectral theorem, but at the same time we can spectrally synthesize this function of $L$ directly from the wave kernel: 
\[
f(\sqrt{L}) = \frac{1}{2\pi} \int_{-\infty}^\infty  \hat{f}(s) \cos (s \sqrt{L})\, ds.
\]
The simple but crucial observation is that this is not just an identity about abstract self-adjoint operators, but also 
calculates the Schwartz kernel of $f(\sqrt{L})$ in terms of the Schwartz kernel of the wave operator.

The following discussion is drawn from the paper \cite{CGT}. Suppose that $f$ has the property that its Fourier transform 
$\hat{f}(s)$ is integrable, along with a certain number of its derivatives, on $\RR \setminus (-\e,\e)$ for any $\e > 0$. 
The first key result is that under such a hypothesis, if $u \in L^2$ has support in a ball $B_r(y)$, then for $R > r$, 
\[
|| f(\sqrt{L}) u||_{L^2(M \setminus B_R(y))}  \leq \pi^{-1} ||u||_{L^2} \int_{R-r}^\infty  |\hat{f}(s)|\, ds. 
\]
The proof is very simple. We know that $\cos (s\sqrt{L}) u$ has support in $B_{r+|s|}(y)$, so that 
\[
||f(\sqrt{L})u|| \leq \frac{1}{\pi} \left\| \int_{R-r}^\infty \hat{f}(s) \cos (s \sqrt{L}) u\, ds\right\| \leq \frac{1}{\pi} ||u|| 
\int_{R-r}^\infty |\hat{f}(s)|\, ds.
\]
A very similar argument gives bounds for $|| L^p f(\sqrt{L}) L^q u||$
depending on the integral of some higher derivatives of $\hat{f}(s)$
over the same half-line. The particularly useful aspect of this is
that the integrals of $|\del_s^\ell \hat{f}|$ which appear on the
right in these estimates start at $R-r$ rather than at $0$, and hence
if these functions decay at some rate, then the right sides of these
inequalities exhibit the corresponding decay.  Assuming we are on a
space with appropriate local uniformity of the metric (or coefficients
of $L$), then we can deduce from this some off-diagonal pointwise
estimates for the Schwartz kernel $f(\sqrt{L})(z,w)$. By off-diagonal
we mean that the estimates are valid in any region where
$\mathrm{dist}\,(z,w) \gg 0$.  One reason for assuming this local
uniformity for $L$ is that these arguments require bounds on the
injectivity radius and volumes of geodesic balls, for example, in
order to pass from $L^2$ to pointwise estimates.

\subsection{Carron's theorem}
This subsection provides a concrete example of how this all works.  We describe some of the main features in the 
paper \cite{Ca} of Carron which uses the ideas above to derive fairly accurate pointwise bounds on the off-diagonal 
decay of the resolvent kernel and heat-kernel for Laplace-type operators on symmetric spaces of noncompact type.

In order to describe this we must first explain at least a small amount about the geometry of these spaces. This 
is recounted elsewhere in much greater detail; the classic reference is \cite{Hel}, but we refer (self-servingly) 
to \cite{MazVas} for an analyst's point of view of this geometry. 

Symmetric spaces are distinguished amongst general Riemannian manifolds by the richness of their isometry
groups.  Their defining property is that the geodesic reflection around any point ($exp_p (v) \mapsto \exp_p(-v)$)
extends to a global isometry; Cartan's classic characterization is that any such space is necessarily of the form 
$G/K$, where $G$ is a semisimple Lie group and $K \subset G$ is a maximal compact subgroup, endowed with an 
invariant metric. Because of this, almost all of the basic structure theory can be reduced to algebra and hence described 
quite explicitly.   We shall focus on one particular realization stemming from the polar decomposition $G = KAK$, 
where $K$ is as above and $A$ is a maximal connected abelian subgroup. For a symmetric space $X$ of noncompact
type, this subgroup $A$ is isomorphic (and isometric) to a copy of $\RR^k$ for some $k$, where the positive integer $k$ 
is called the rank of $X$. Using this polar decomposition, we identify $G/K \cong KA$. The map $\Phi: K \times A \to X$, 
$\Phi(k,a) = ka$, is surjective, but far from injective. 

It is best to think of the simplest special case, the real hyperbolic space $\HH^n$; here $A \cong \RR$ and 
$K = \mbox{SO($n$)}$. The image of the origin $0 \in A$ via $\Phi$ is a single point $o \in X$, and this point
is fixed by the entire (left) action of $K$.  The space $X$ is the union of geodesic lines through $o$ which all intersect 
pairwise only at this point. The group $K$ acts transitively on this space of geodesic lines through $o$ with stabilizer 
$\mbox{SO($n-1$)}$.  Note that there are elements of $K$ which take a geodesic to itself but reverses its orientation;
this means that we get a less redundant `parametrization' by restricting $\Phi$ to $K \times \RR^+$.
Geometrically, we have the familiar picture of $\HH^n$ as $\RR^+ \times S^{n-1}$ with the warped product metric 
$dr^2 + \sinh^2 r \, d\theta^2$. 

For a general symmetric space $X$ of rank $k > 1$, this picture generalizes as follows. The space $X$ is the 
union of the various images of $A$ by elements $k \in K$, and all of these images intersect at $o$, though
$kA \cap k' A$ often consists of a larger subspace. These translates of $A$ by $k \in K$ should be thought
of as the radial directions in $X$.   Another important piece of structure is the existence of a finite set of linear 
functionals $\Lambda = \{\alpha_j\}$ on $A$ called the roots. These divide into positive and negative roots,
$\Lambda = \Lambda^+ \cup \Lambda^-$, and the positive roots determine a (closed) sector $V \subset A$ by
$V = \{\alpha_j \geq 0 \ \forall\  \alpha_j \in \Lambda^+\}$. This sector $V$ is the analogue of the half-line
in $\HH^n$, and the restriction of $\Phi$ to $K \times V$ is still surjective, and if $K' \subset K$ denotes the isotropy
group at a generic point, then we can regard $X$ as being the product $V \times (K/K')$ with certain
submanifolds of $K/K'$ collapsed along various boundary faces of $V$. In terms of this data, we can finally write 
down the multiply-warped product metric  
\[
g = da^2 + \sum_j \sinh^2 \alpha_j \, dn_j^2,
\]
where the sum is over positive roots, $da^2$ is the Euclidean metric on $A$ and $dn_j^2$ is a metric on a certain 
subbundle of the tangent bundle of $K/K'$ corresponding to the root $\alpha_j$. 

For simplicity here we just discuss the scalar Laplacian $\Delta$ on $X$ and, following the theme
of this section, consider the problem of estimating the Schwartz kernels of $f(\sqrt{-\Delta})$ for
suitable functions $f$.  Because $\Delta$ commutes with all isometries on $X$, the Schwartz kernel 
$K_f(x, x')$ of this operator depends effectively on a smaller number of variables. Given any pair
of distinct points $x, x' \in X$, choose an isometry $\varphi$ of $X$ so that $\varphi(x') = o$ and $\varphi(x)$ 
lies in some particular copy of $A$. If we ask that $\varphi(x) = a \in V \subset A$, then $\varphi$ is 
almost uniquely determined.  We thus have that $K_f(x, x') = K_f(\varphi(x), o) = K_f(a,0)$. In other words, 
$K_f$ is really only a function of the $k$ Euclidean variables $a = (a_1, \ldots, a_k)$. In particular, when 
the rank of $X$ is $1$, then $K_f$ reduces to a function of one variable $r \geq 0$.

This reduction points to the difficulty of studying functions of the Laplacian on symmetric spaces of
rank greater than $1$. Indeed, while the resolvent kernel $R(\lambda; x, x')$ on a space of rank $1$
depends only on $\mbox{dist}(x,x')$ and hence can be analyzed completely by ODE methods, the
same is not true when the rank of $X$ is larger. Similarly, even in the rank $1$ setting, the heat
kernel $H(t, x, x')$ depends on two variables, $t$ and $\mbox{dist}(x,x')$, but unlike the Euclidean
case, there is no extra homogeneity which reduces this further to a function of one variable. Thus
the problem which Carron's theorem answers is how to give good estimates on these reduced functions,
$R(\lambda, a)$ and $H(t, a)$, where $a \in A$ is the `relative position' between $x, x' \in X$. 

\begin{theorem}[Carron \cite{Ca}]
Let $X$ be a symmetric space of noncompact type and rank $k$ and consider the Schwartz kernels 
$R(\lambda, a)$ and $H(t, a)$ of the resolvent $(-\Delta - \lambda_0 - \lambda^2)^{-1}$ and heat operator $e^{t\Delta}$,
written in reduced form as above. The number $\lambda_0$ here is the bottom of the spectrum of $-\Delta$;
it may be calculated explicitly.  Then
\[
|R(\lambda, a)| \leq C e^{-\rho(a) - \mathrm{Re}(\lambda)\, \operatorname{dist}(a,o)}
\]
and
\[
|H(t,a)| \leq C e^{-\lambda_0 t - \rho(a) - \operatorname{dist}^2(a,0)/4t} \, \Phi_t(a).
\]
The function $\Phi_t(a)$ is a somewhat messy but quite explicit and understandable function which is 
a rational function of $a$ and certain powers of $t$. The linear functional $\rho$ on $A$ is half the sum of 
the `restricted' positive roots; this is a standard object which appears frequently. 
\end{theorem}
It is known that the upper bounds given here are sharp in the sense that there are lower bounds that
differ just by the constant multiple for these same kernels. We refer also to the papers \cite{AJ} and \cite{LM}
sharper bounds obtained by different and more complicated methods.

The proofs of these estimates are clever but not very long, and in the remainder of this section we give a few of the 
ideas which go into them. 

The first step is that if $a \in A$ is arbitrary and $\e \in (0,1)$, then we can estimate from above and below the 
volume of the set $K B(a,\e)$, where $B(a,\e)$ is a ball of radius $\e$ in $A$ centered around $a$. The reason this
can be done is because we have very good information on the Jacobian determinant for the coordinate change
implicit in some natural coordinatizations induced by $K \times A \to X$. 

Let us first study the resolvent. Fix $a \in A$ such that $\mbox{dist}(a,o) \geq 2$. Let $D$ be the ``annular''
region $K B(a,1)$ and
suppose that $\sigma \in L^2(D)$ (so $\sigma$ vanishes outside $D$).  We estimate $R(\lambda, a)\sigma$ 
and then $R(\lambda, a)$ itself, outside this region, i.e. where $\mbox{dist}(a,o) \geq 2$. If $u = R(\lambda, \cdot) \sigma$, 
then 
\[
u = \int_0^\infty \frac{e^{-\lambda \xi}}{\lambda} \cos (\xi \sqrt{-\Delta - \lambda_0}) \sigma \, d\xi.
\]
Then using that $||\cos (\xi \sqrt{-\Delta - \lambda_0})||_{L^2 \to L^2} = 1$ as well as finite propagation
speed, because of the support properties of $\sigma$, we obtain that
\[
||u||_{L^2(B(o,1))} \leq \int_{\operatorname{dist}(a,o) - 2}^\infty 
\frac{e^{ -\mathrm{Re}(\lambda) \xi}}{|\lambda|} ||\sigma||_{L^2} \leq 
\frac{1}{|\lambda|^2} e^{-\mathrm{Re}(\lambda) (\operatorname{dist}(a,o) - 2)} ||\sigma||_{L^2}.
\]
From here, using local elliptic estimates, we obtain that
\[
|u(o)| \leq C e^{-\mathrm{Re}(\lambda) \operatorname{dist}(a,o)} ||\sigma||_{L^2(D)}.
\]
In other words, this estimates the norm of the mapping $T$ defined by $L^2(D) \ni \sigma \mapsto 
\left. R(\lambda)\sigma\right|_{o}$, whence (using the $L^\infty \to L^\infty$ norm of $TT^*$),  
\begin{equation}
\int_{D} |R(\lambda, x, o)|^2 \, dx \leq C e^{-2\mathrm{Re}(\lambda) \operatorname{dist}(a,o)}. 
\label{intestR}
\end{equation}

We next wish to find a similar estimate where the integral on the left is only over some ball $B(ka, 1/4)\subset D$
rather than the entire annular region $D$.  More specifically we assert that
\[
\mbox{Vol}(B(ka,1/4)) \, \int_{B(ka, 1/4)} |R(\lambda, x, o)|^2\, dx \leq C e^{-2\mathrm{Re}(\lambda) \operatorname{dist}(a,o)}. 
\]
This must be true since if it failed for all such balls then the sum over all balls would lead to a violation of \eqref{intestR}.

Finally, noting that the volume of this ball is approximately $e^{2\rho(a)}$, and applying the same local
elliptic estimates as before to estimate the value at a point in terms of a local $L^2$ norm, we conclude that
\[
|R(\lambda, a, o)| \leq C e^{-\mathrm{Re}(\lambda) \operatorname{dist}(a,o) - \rho(a)}. 
\]
This is the desired off-diagonal decay estimate. 

The corresponding argument to estimate the off-diagonal behaviour of the heat kernel proceeds in a very
similar way, substituting local parabolic estimates for local elliptic estimates. We refer to \cite{Ca} for details.

It is worth remarking that there are other very effective ways to
establish so-called Gaussian bounds for heat kernels under rather
general circumstances. We mention in particular the beautiful theory
developed by Grigor'yan and Saloff-Coste, see \cite{Grig}, \cite{SC}.
These techniques work in far more general circumstances, and depend on
quite different underlying principles. However, one point of interest
in Carron's work is that he is able to obtain the correct polynomial
factors in front of the exponential, which is perhaps impossible using
those more general approaches.

\section{Scattering theory}
For the second and longer part of this survey, we turn to an entirely different aspect of the global theory of wave equations 
and discuss some approaches to mathematical scattering theory.  This classical subject has deep physical origins, and
has received numerous mathematical formulations. While these approaches are mostly equivalent, the correspondences 
between them are sometimes not always obvious.  In the following pages we first review one point of view on stationary 
scattering theory, then turn to some perspectives on the corresponding time-dependent theory.  This is all done with
a distinctly PDE (rather than, say, operator-theoretic) focus. We conclude with a discussion of a more abstract functional 
analytic setup of scattering theory due to Lax and Phillips centered around the notion of a translation representation and
explain how the theory of radiation fields developed by Friedlander provides a concrete realization of the translation
representation. 

There are numerous settings in which to introduce any of these topics, including scattering by potentials, which is 
the study of Schr\"odinger operators $-\Delta + V$ on $\RR^n$, or scattering by obstacles, which studies these same
operators but on exterior domains $\RR^n \setminus \mathcal O$ with some elliptic boundary condition at $\partial \mathcal O$. 
There are also significant differences between these theories when the ambient space is odd- or even-dimensional.  
Finally, it is also natural to consider these same problems on manifolds which are asymptotically Euclidean or 
asymptotically conic at infinity (or indeed, have some other type of asymptotically regular geometry, e.g.\ 
asymptotically hyperbolic). Each setting requires different sets of techniques, and in order to make this 
exposition as simple as possible, we focus on the combination of hypotheses where everything works out
most simply. Namely, we study the scattering theory associated to $L = -\Delta + V$ on an {\it odd}-dimensional 
Euclidean space $\RR^n$, with the strong assumption that $V \in \calC^\infty_0$.  We describe the structure 
theory for solutions of the Helmholtz equation $(L - \lambda^2)u = 0$, 
and for $\Box_V := \Box + V = D_t^2 - L$, the time-dependent wave equation, and give some indication
how objects in these respective settings correspond to one another.

There are very many excellent references to each part of what we discuss (and much that is closely related that
we do not discuss), so we relegate almost all of the technicalities to those sources. We mention in particular  
\cite[Vol. IV]{RS}, \cite[Ch. 9]{Taylor}, \cite{Perry-Enss}, \cite{Yafaev} and \cite{Melrose}. The material on 
radiation fields is spread over several papers, starting from the original work by F.G.\ Friedlander
\cite{Fr}.  There is a forthcoming and detailed survey of this subject by Melrose and Wang \cite{MW}, to which 
the discussion here is intended to be an introduction. 

\subsection{Stationary scattering theory}
The stationary formulation of scattering theory concerns the elliptic operator $L - \lambda^{2}$, where
here and below, $L = -\Delta + V$, with $V \in \calC^\infty_0$ (and real-valued!). 
It is obvious that $L$ is bounded below, i.e.\ 
\[
\int_{\RR^n} (Lu) \overline{u}\, dV \geq -C \int_{\RR^n} |u|^2\, dV
\]
for all $u \in \calC^\infty_0(\RR^n)$, and with little more work one can also prove that it has a unique self-adjoint 
extension as an unbounded operator on $L^2(\RR^n)$.  Indeed, this is yet another consequence of the 
finite speed of propagation, see \cite{Chernoff}. Its spectrum is contained in a half-line $[-C, \infty)$; 
the positive ray $[0,\infty)$ comprises the entire continuous spectrum, and there are a finite number of $L^2$ 
eigenvalues in the $[-C, 0)$. If we allow $V$ to be less regular, simple examples show that this negative 
interval may contain an infinite sequence of such eigenvalues converging to $0$; the basic example of
this is when $V(x) = -1/|x|$, which is the potential for the Schr\"odinger operator modeling the hydrogen atom. 

Assume initially that $\lambda$ lies in the lower half-plane $\Im \lambda < 0$. Provided that $\lambda \neq -i |\lambda_j|$
corresponding to any of the negative eigenvalues $-\lambda_j^2 < 0$, the operator $L -\lambda^2$ has an $L^2$ 
bounded inverse, 
\[
R_V(\lambda) = (L - \lambda^2)^{-1}.
\] 
This is called the resolvent and is a meromorphic family of bounded operators on $L^2$ with simple poles at 
the points $-i|\lambda_j|$. 

The first issue is to show that the continuous spectrum
($\lambda^{2}\in [0,\infty)$) is absolutely continuous, or in other words,
that the singular continuous part of the spectrum is empty. More specifically, we must find an $L$-invariant 
orthogonal splitting $L^2(\RR^n) = \calH_{\mathrm{pp}} \oplus \calH_{\mathrm{ac}}$, so that the restriction of $L$ 
to $\calH_{\mathrm{pp}}$ is discrete, while the restriction of $L$ to ${\calH}_{\mathrm{ac}}$ is absolutely continuous. 
It is a classical theorem due to Friedrichs that in this setting any $L^2$ eigenvalue of $L$ is strictly negative.
The proof consists of showing that if any such eigenvalue is positive, then the corresponding eigenfunction 
must vanish outside a compact set, which violates standard unique continuation theorems. (This uses that 
$V$ is compactly supported -- if $V$ only decays rapidly then the argument is a bit more intricate.) 
By the general spectral theorem, the absolute continuity of $\left. L  \right|_{\calH_{\mathrm{ac}}}$ is equivalent
to the existence of a unitary isomorphism $U: \calH_{\mathrm{ac}} \longrightarrow L^2(\RR;Y)$, where $Y$ is an 
auxiliary Hilbert space, so that the self-adjoint operator $U \circ L \circ U^{-1}$ on $L^2(\RR; Y)$ is
multiplication by the coordinate function $t \in \RR$. One of the goals of scattering theory is to 
exhibit this unitary isomorphism explicitly, which is done using the M{\o}ller wave operators, see below. 
A closely related goal is to understand the structure of generalized (non-$L^2$) solutions to the equation 
$(L - \lambda^2) u = 0$, $\lambda^2 > 0$.  The key tool for all these questions is the resolvent $R_V(\lambda)$, introduced above. 



Let us first consider the free Laplacian $L_0 = -\lap$ on $\reals^{n}$. When $\lambda \in \reals \setminus\{0\}$, 
the nullspace $\mathcal E(\lambda)$ of the operator $-\lap - \lambda^{2}$ contains the plane wave solutions
$e^{i\lambda z   \cdot \omega}$ for any $\omega \in \sphere^{n-1}$. Any linear combination of these plane waves 
also lies in $\calE(\lambda)$, and indeed, general superpositions of these plane wave solutions span all of
$\calE(\lambda)$.  We explain this more carefully. For any $g \in \calC^{\infty}(\sphere^{n-1})$, define 
\begin{equation*}
u( z) = \int _{\sphere^{n-1}}e^{i\lambda z \cdot \omega}g(\omega) \,d\omega.
\end{equation*}
This is a solution of $(-\lap - \lambda^{2})u = 0$, and the most general (polynomially bounded) element of 
$\calE(\lambda)$ can always be obtained from this same representation but allowing $g$ to be a distribution. 
The ``smooth'' elements of $\calE(\lambda)$ are those where $g$ is smooth. 

We can look at this a different way. Note that since $\omega \mapsto z \cdot \omega$ is a Morse function 
on $\sphere^{n-1}$, and has critical points $\omega = \pm z / |z|$, the stationary phase lemma
shows that (assuming $g$ is smooth), the integral expression for $u$ has an asymptotic expansion of the form 
\begin{equation}
\label{eq:expansion}
u(z) \sim e^{i\lambda|z|} |z|^{-\frac{n-1}{2}}\sum _{j  =0}^{\infty}|z|^{-j}a_{+,j}(\theta) +
e^{-i\lambda|z|}|z|^{-\frac{n-1}{2}}\sum_{j=0}^{\infty}|z|^{-j} a_{-,j}(\theta).
\end{equation}
Here $z = |z|\theta$, $\theta \in \sphere^{n-1}$ are polar coordinates on $\RR^n$.  As part of this, one obtains 
that up to a multiple of $2\pi$, $a_{\pm,0}  = i ^{\mp  (n-1)/2}g(\pm \theta)$.   Closely related is the assertion that
any $u \in \calE(\lambda)$ has an expansion of this same form and moreover, fixing any $a_{+,0}\in \mathcal D'(S^{n-1})$,
there is a unique $u \in \calE(\lambda)$ with this distribution as its leading coefficient.
It is reasonable to regard the operator $\calP: a_{+,0} \mapsto u$ as solving a Dirichlet problem at infinity for 
$-\Delta - \lambda^2$, and hence we call $\calP$ the \emph{Poisson operator}. 

The free \emph{scattering operator} at energy $\lambda$ is the map $\mathcal{S}_0(\lambda)$ sending the function 
$a_{+,0}$ to $a_{-,0}$.  Using the explicit representation above, we see that in this free setting, 
$\mathcal{S}_0(\lambda) a (\theta) = i^{n-1}a(-\theta)$; it is just a constant multiple of the antipodal map.

Proceeding slightly further with the free problem, suppose that $\im \lambda < 0$. Using the Fourier transform,
one can determine the inverse of $-\lap - \lambda^{2}$ (as an operator on Schwartz functions) via
\begin{equation*}
R_{0}(\lambda) f = (-\lap - \lambda^{2} )^{-1} f = (2\pi)^{-n}\int _{\reals^{n}} e^{iz\cdot \zeta}(|\zeta|^{2} - \lambda^{2} )^{-1}
\hat{f}(\zeta) \,d\zeta.
\end{equation*}
When $n$ is odd, this has a particularly simple form: there is a simple polynomial $p_{n}(\alpha)$ of degree $(n-1)/2$ such that 
the Schwartz kernel of $R_{0}(\lambda)$ can be written as 
\begin{equation}
\label{eq:resolvent-odd-dim}
|z-z'|^{2-n}\, p_{n}(\lambda|z-z'|) e^{-i\lambda|z-z'|}.
\end{equation}
(In particular, $p_3(\alpha)$ is simply a constant.) 
There is a related but slightly more complicated formula when $n$ is even. This explicit expression shows that 
as a function of $\lambda$, $R_{0}(\lambda)$ continues holomorphically from the lower ``physical'' half-plane 
$\{ \Im \lambda < 0 \}$ to the entire complex plane when $n\geq 3$ is odd. When $n=1$, this continuation has a simple pole at $\lambda = 0$, 
and when $n$ is even, there is a similar continuation but to the infinitely sheeted logarithmic Riemann surface 
branched at the origin.  To make sense of this, one can say that this Schwartz kernel continues as a holomorphic 
function taking values in distributions;  an alternate and equivalent sense is to regard the continuation 
taking values in the space of bounded operators $L^{2}_{c} \to L^2_{\mathrm{loc}}$, (this domain space consists 
of compactly supported $L^{2}$ functions). From \eqref{eq:expansion} and stationary phase, one proves that if $f\in
\calC^{\infty}_{c}(\reals^{n})$, then 
\begin{equation*}
R_{0}(\lambda) f = e^{-i\lambda|z|} |z|^{-\frac{n-1}{2}}w,
\end{equation*}
where $w$ is a smooth function on the radial compactification of $\reals^{n}$. This last assertion about smoothness
on the compactification is simply a concise way of stating that $w$ has an asymptotic expansion
\[
w \sim \sum_{j=0}^\infty w_j(\theta) |z|^{-j}.
\]

Let us now pass to the analogous considerations for the operator $L$. Some versions of all the structural results about 
solutions remain true. These are typically proved by a perturbative argument, which means that one no longer 
has explicit formul\ae.  The starting point is the Lipman--Schwinger formula, which gives a relationship between
$R_0(\lambda)$ and $R_V(\lambda)$ in the region in the $\lambda$-plane where they both make sense. This states that
\begin{equation*}
R_V(\lambda) = R_0(\lambda) \left( I + VR_{0}(\lambda)\right)^{-1} = 
\left( I + R_{0}(\lambda) V\right)^{-1} R_0(\lambda). 
\end{equation*}
The issue is to prove that the inverses of $I + V R_0(\lambda)$ and $I + R_0(\lambda)V$ make sense, and to do
this one observes that $VR_{0}(\lambda)$ and $R_0(\lambda) V$ are compact operators (between suitable
function spaces), so that one can invoke the analytic Fredholm theorem to obtain that these inverses,
and hence $R_V(\lambda)$ itself, are meromorphic on the region where $R_0(\lambda)$ is holomorphic (hence
on $\CC$ when $n$ is odd and greater than $1$). 

The argument sketched earlier that $L$ has no $L^2$ eigenvalues embedded in the continuous spectrum implies 
that $R_V(\lambda)$ has no poles on the real axis.  (The argument for regularity at $\lambda = 0$ requires slightly 
more care.)  On the other hand, the negative eigenvalues $\lambda_j$ of $L$ correspond to poles of $R_V(\lambda)$
at $-i |\lambda_j|$.  The new and perhaps unexpected phenomenon in that $R_V(\lambda)$ may have poles in 
the upper half-plane (and indeed, this always occurs if $V$ is nontrivial).  These poles are known as the resonances 
of $L$, and their location and statistical distribution has been the target of much research. 


Let $\calE_V(\lambda)$ denote the nullspace of $L - \lambda^{2}$ (say in $\calS'(\RR^n)$).  Just as in the free case,
this space may be generated using ``distorted'' plane waves. These are defined as follows. For any 
$\omega \in \sphere^{n-1}$ and $\lambda \in \reals \setminus\{0\}$, there is a function $W_{\lambda, \omega}$ which is
smooth on the radial compactification of $\reals^{n}$ so that 
\begin{equation*}
\phi_{\lambda,\omega}(z) = e^{i\lambda z\cdot \omega} + e^{-i\lambda|z|}|z|^{-\frac{n-1}{2}} W_{\lambda, \omega} 
\end{equation*}
lies in $\calE_V(\lambda)$. Note that the second term here is simply $R_{0}(\lambda)(-V e^{i\lambda z \cdot  \omega})$.
Superpositions of these can be used as before to generate all elements of $\calE_V(\lambda)$. Indeed, 
if $g\in \calC^{\infty}(\sphere^{n-1})$, then the general ``smooth'' element of $\calE_V(\lambda)$ can be written as
\begin{equation*}
u(z) = \int _{\sphere^{n-1}} \phi_{\lambda,\omega}g(\omega)\,d\omega,
\end{equation*}
Using stationary phase as before, this integral has an asymptotic expansion of exactly the same form as 
\eqref{eq:expansion}. The leading coefficient $a_{+,0}(\theta)$ is again just (a multiple of) $g$, but now
the other leading coefficient $a_{-,0}(\theta)$ is not simply the reflection $g(-\theta)$, but rather
a sum of this reflection plus an extra term which is an integral over $S^n$ involving both $g$ and $V$. 
The scattering operator $\calS_V(\lambda)$, which sends $a_{+,0} \mapsto a_{-,0}$, is again unitary, and
is the sum of the antipodal operator and another term which has a smooth Schwartz kernel. 
The map $\calP_V(\lambda)$ which sends $a_{+,0}$ to $u$ is again called the Poisson operator.

The results and definitions above continue to hold in suitably modified form not only for obstacle scattering, 
but also in the rather general setting of asymptotically Euclidean or asymptotically conic manifolds 
(these are called \emph{scattering manifolds}~\cite{Melrose:1994} by Melrose).  For more on this as well as 
many further details about everything discussed above, we refer to the book of Melrose~\cite{Melrose},
see also \cite{Melrose:1994} and \cite{MZ}.  

\subsection{Time-dependent scattering}
\label{sec:time-depend-scatt}
We now turn our attention to the time-dependent formulation of scattering theory, 
and its relationship with stationary scattering.  This time-dependent theory involves
the study of ``large time'' properties of solutions of the wave equation.  The connection
with the stationary approach is via the Fourier transform in time; indeed, this
Fourier transform carries $L - D_t^2$ to $L - \lambda^2$, and asymptotic properties
of as $|t| \to \infty$ correspond to `local in $\lambda$' properties of the latter operator.
For the wave equation associated to $L = -\Delta + V$, where $V$ is compactly supported,
the intuitive picture is that one sends in a wave for times $t \ll 0$ from some direction 
at infinity and then observes what happens as this wave interacts with the potential 
and then scatters into a sum of plane waves as $t \nearrow +\infty$.  Amongst the many
good sources for this material, we refer to the books of Friedlander~\cite{Friedlander:wave}, 
Lax~\cite{Lax}, Lax--Phillips~\cite{LP}, Taylor~\cite{Taylor} and Melrose~\cite{Melrose}.

\subsubsection{Progressing wave solutions}
\label{sec:progr-wave-solut}
We begin by describing the special class of progressing wave solutions for wave operators.
The calculations here go back to the dawn of microlocal analysis and can be regarded as
the nexus of many constructions and ideas in that field. 
This construction is quite geometric and it is most naturally phrased in terms of the wave operator
on a general Lorentzian metric $g$. The special case of a static metric $g = -\dt^{2} + h$ 
on the product of $\RR$ with a Riemannian manifold $(M,h)$ is of particular interest, and
we discuss at the end how this specializes for the particular operator $\Box_V = \Box + V$ on Minkowski space.

Thus let $(X,g)$ be a Lorentzian manifold and consider $\Box _{g} + V$, where $V \in \calC^\infty_c(X)$. 
We look for solutions $u$ to $(\Box_g + V)u = 0$ which have the form
\begin{equation*}
u = \coeff \, \wf (\ph),
\end{equation*}
where $\coeff$ is smooth, $\wf$ is a distribution on $\reals$ which models the `wave form' of the solution, 
and $\ph$ is a function on $X$ with nowhere vanishing gradient which we call the phase function.
To be concrete, we typically let $\wf = \delta$ or $\wf = x^{k}_{+}$ for some $k \geq 0$, but the key feature
we require of $\wf$ is that it behave like a homogeneous function in the sense that its successive derivatives and
integrals are progressively more or less smooth than $\wf$ itself.  Of course, it is usually impossible to choose 
solutions of $(\Box_g + V)u = 0$ which have this precise form, but the goal is to add increasingly higher order 
correction terms of a similar form involving the integrals of $\wf$ so that, in the end,
this initial expression is the first term in some asymptotic expansion of an exact solution. 

The first step is to calculate
\begin{multline*}
\left( \Box_{g} + V\right) u = \frac{1}{\sqrt{|g|}}\pd[i]\left(g^{ij}\sqrt{|g|}\pd[j]u\right) + Vu \\
= \wf ''(\ph) g\left(\grad \ph, \grad \ph\right)\coeff + \wf'(\ph) \left( 2g\left(\grad\ph, \grad \coeff\right)
+ \coeff \Box_{g}\ph\right) \\ + \wf (\ph) \left( \Box \coeff + V\coeff\right)
\end{multline*}

As indicated above, assume that $\wf_{k}$ is a sequence of distrubutions on $\RR$ such 
that $\wf_{k} =\wf_{k+1}'$.  (Again, refer to the basic example $\wf_{0} = \delta$, $\wf _{k+1} = 
\frac{1}{k!}x_{+}^{k}$.)  Let us now assume that 
\begin{equation}
u \sim \sum_{k\geq 0} u_{k} = \sum _{k\geq 0}\coeff_{k}(t,z) \wf_{k}(\ph).
\label{sumseries}
\end{equation}
We apply the calculation above and group together the terms of the same order (where the order
of $\alpha_k$ is $k$ and each derivative lowers the order by $1$).  

Grouping terms of the same order, we attempt to choose $\coeff_{k}$ so that each term vanishes.  The only term of order $-2$ is 
$\coeff_{0}\wf_{0}''(\ph) g\left( \grad\ph,  \grad\ph\right)$, so the first requirement is that
\begin{equation*}
g \left( \grad\ph, \grad \ph\right) = 0.
\end{equation*}
This is known as the \emph{eikonal equation} and states that $\grad \ph$ is a \emph{null-vector} for
the metric $g$.  This is a global nonlinear Hamilton-Jacobi equation for $\ph$.  In the special case 
$X = \reals \times M$, $g = -\dt^{2} + h$, the eikonal equation can be written as 
\begin{equation*}
  \left( \pd[t]\ph\right)^{2} = \left| \grad_{h}\ph\right|^{2};
\end{equation*}
if we write $\ph = t - \phr$, where $\phr$ is a function on $M$, then 
\begin{equation*}
  \left| \grad_{h} \phr\right|^{2}= 1.
\end{equation*}
It is straightforward to see that the level sets $\phr = \mbox{const}$ are 
at constant distance from one another, so in general, $\phr(x) = \mbox{dist}_h\, (x, Z)$
where $Z$ is some fixed level set of $\phr$.  Even in the more general Lorentzian
setting, the function $\ph$ incorporates a lot of the distance geometry of $g$. 

In any case, fix a solution $\ph$ of the eikonal equation. We have now arranged 
that the term of order $-2$ vanishes. In fact, for any $k$, the term
containing $g(\grad\ph, \grad \ph)$ vanishes, and so the
equations for the higher coefficients simplify to \emph{transport  equations}.  
In particular, the term of order $-1$ reduces to 
\begin{equation*}
  \wf_{0}'(\ph) \left( 2 g \left( \grad \ph, \grad \coeff_{0}\right) +
    \coeff_{0}\Box \ph\right).
\end{equation*}
Since $\grad\ph$ is nowhere vanishing, this is a linear ODE for $\coeff_{0}$ along the integral curves 
of $\grad \ph$, which means that given any initial conditions for $\coeff_{0}$ on the 
characteristic surface $\ph = \text{constant}$ we may solve this equation locally. 

The term of order $k-1$ yields an inhomogeneous transport equation
for $\coeff_{k}$ in terms of $\ph, \coeff_{0},\ldots, \coeff_{k}$.  We
solve this transport equation with vanishing initial data and proceed
inductively to choose all $\coeff_{k}$. 


It is possible to asymptotically sum the series \eqref{sumseries}. This means that we can
choose a function $v$ with the property that
\[
v - \sum_{k=0}^N \coeff_k \wf_k(\ph)
\]
is as smooth as the next term in the series, $\coeff_{N+1} \wf_{N+1}(\ph)$.  By construction,
$(\Box_g + V) v = f \in \calC^\infty(X)$.  We must now invoke a theorem guaranteeing the
existence of a smooth solution $w$ for the initial value problem $(\Box_g + V) w = f$ with
vanishing Cauchy data vanishes, where $f$ is smooth. Given this, then $u = v - w$ 
is a solution of the original equation and the expansion we have calculated determines
the singularity profile of $u$. Note that these singularities of $u$ occur precisely along
the union of level sets $\ph = c$ where one (and hence every) $\alpha_k$ is singular
at $c$. 

For the special case where $g = -dt^2 + dx^2$ on Minkowski space, fix $\omega \in \sphere^{n-1}$ and 
consider the equation 
\begin{equation*}
 \label{forcing}
\left(\pd[t]^{2} - \lap_{z} + V\right) u = 0, \qquad   u = \delta (t - z\cdot \omega)\ \ \mbox{when}\ t \ll 0. 
\end{equation*}
The eikonal equation $|\grad \ph|_g^{2} = 0$ has solution $\ph(t,z) = t - z\cdot \omega$.  This gives a 
global solution of the wave equation for all $t$ when $V \equiv 0$. However, by the propagation of singularities 
theorem, the wave front set of the solution $u$ for the perturbed problem with this initial data in the distant past
agrees with that of this exact free solution. Hence it makes sense to look for a solution of the perturbed
problem of the form 
\begin{equation*}
u \sim \delta (t - z\cdot \omega) + \sum _{k\geq    0}\coeff_{k}(t,z)  x_{+}^{k}(t-z\cdot \omega),
\end{equation*}
for some choice of smooth functions $\coeff_{k}$.  This fits exactly into the scheme above (and was,
of course, the setting for the original version of these calculations). The first transport equation is 
\begin{equation*}
2\left( \pd[t] - \omega \cdot \grad_{z}\right)\coeff_{0} = 0,
\end{equation*}
which means that $\coeff_0$ is a function of $t = z \cdot \omega$ and $z$; its Cauchy data is 
defined on the hypersurface $t = z\cdot \omega$, and the equation dictates
that it must be constant along the lines parallel to $\omega$. 

Once we have determined $\coeff_{0}, \ldots , \coeff_{k}$, then the $(k+1)^{\mathrm{st}}$ transport equation is 
\begin{equation*}
 2(k+1)\left( \pd[t] - \omega \cdot \grad_{z}\right) \coeff_{k+1} =   -\left( \Box  + V\right) \coeff_{k},
\end{equation*}
which we solve with vanishing initial data.  Carrying this procedure out for all $k$ determines the Taylor series of $u$ 
along the hypersurface $\{t = z \cdot \omega\}$. As described earlier, we can use the Borel Lemma to choose
an asymptotic sum $v$ for this series, so that $\left( \Box + V\right) v = f$ is smooth and $v$ satisfies the
correct ``initial condition'' for $t \ll 0$. We can then find a \emph{smooth} correction term $w$ which solves
away this error term. Thus $u = v - w$ is an exact solution

The calculations here were historical precursors to the more elaborate but ultimately very similar ones which 
come up in the construction of Fourier integral operators.  Indeed, solving the eikonal equation for $\ph$ 
is the direct analogue of solving the eikonal equation for the phase of an FIO. For potential scattering, 
keeping track of the parametric dependence on $\omega$ fixes the phase; the solutions of the transport 
equations are the coefficients in the expansion of the amplitude, and these correspond to 
the terms in the expansion for the symbol of the FIO.

\subsubsection{M{\o}ller wave operators}
\label{sec:moll-wave-oper}
We now turn to another perspective on time-dependent scattering, which is through the definition
of the so-called M{\o}ller wave operators. This can be regarded as a formalization of the discussion above;
there we described how to calculate the profile of the solution obtained by ``sending in''  a delta function 
along a particular direction. Our goal now is to put this information together into a map which compares
the long-time evolution with respect to the perturbed equation against that for the free equation.

Let us suppose now that $g = -dt^2 + h$ is a static Lorentzian metric. For any ($\calC^\infty_c$) potential $V$, 
define the wave evolution operator
\begin{equation*}
U_{V}(t): C^{\infty}_{c}(\reals^{n}) \times C^{\infty}_{c}(\reals^{n})\to C^{\infty}_{c}(\reals^{n}) \times C^{\infty}_{c}(\reals^{n}),
\end{equation*}
where, if $u$ solves the Cauchy problem 
\[
 \left( \Box + V \right) u =0, \quad  \left( u , \pd[t]u \right) |_{t=0} = (\phi, \psi),
\]
then
$U_V(t_0)(\phi, \psi) = (u, \pd[t]u)|_{t=t_0}$.  The free wave evolution operator $U_{0}(t)$ is defined analogously
using solutions for $\Box u = 0$ instead.  Uniqueness of solutions of these Cauchy problems implies that
$U_{V}$ and $U_{0}$ are groups, i.e.\ $U_*(t)^{-1} = U_*(-t)$ and $U_*(t+s) = U_*(t)U_*(s)$ for $* = 0$ or $V$. 

Now define the \emph{M{\o}ller wave operators} $W_{\pm}$ by 
\begin{equation*}
  W_{\pm}(\phi, \psi) = \lim_{t\to \pm\infty} U_{V}(-t)   U_{0}(t)(\phi, \psi), 
\end{equation*}
presuming the limit exists. This limit is meant to be taken in the sense of strong operator
convergence.  If we define the \emph{energy space}
\[
H_E = \left\{(\phi, \psi): \int \psi^2 + |\nabla_z \phi|^2\, dV  < \infty \right\},
\]
then $W_{\pm}$ extends by continuity to all of $H_{E}$. It can be
proved that if certain local measurements of this energy decay
appropriately, then $-\lap + V$ has no $L^{2}$ eigenvalues and this
extension is an isomorphism of $H_{E}$ to itself.  If $-\lap + V$ does
have $L^{2}$ eigenvalues, then $\calH_{\mathrm{pp}}$ determines a
finite dimensional subspace in $H_E$ and $W_{\pm}$ is an isomorphism
from $H_{E}$ onto the orthogonal complement of $\calH_{\mathrm{pp}}$,
which we denote $H_E^\perp$.

Since $U_0(t)$ and $U_V(t)$ are unitary, the wave operators $W_{\pm}$ are characterized by the property 
that
\[
\norm[H_{E}]{U_{V}(t)W_{\pm}(\phi, \psi) - U_{0}(t)(\phi, \psi)}
\to 0 \text{ as }t\to \pm \infty.
\]
for all $(\phi,\psi) \in H_E^\perp$.  Now define the \emph{scattering operator} 
\begin{equation*}
\mathcal{S} = W_{+}^{-1}W_{-}; 
\end{equation*}
this is an isomorphism of $H_{E}^\perp$.  It describes the relationship between the 
asymptotic free wave emerging as $t \nearrow +\infty$ for a solution of the perturbed
equation $(\Box + V)u = 0$ in terms of the incoming free wave for $t \ll 0$. 

These operators lead directly to the unitary isomorphism mentioned earlier which intertwines $L$ (or rather,
its restriction to $\calH_{\mathrm{ac}}$), with a simple multiplication operator. In other words, the existence and
properties of the wave operators and scattering matrix proves that the singular continuous spectrum of $L$ is empty.

There are many other settings where one can define analogues of the  M{\o}ller wave and scattering operators.  
Classically this is done for exterior domains, and more recently on asymptotically Euclidean or conic
manifolds (where the structure of the scattering matrix is quite intriguing, see \cite{MZ}), as well as other
geometric settings such as asymptotically hyperbolic manifolds, etc.  There is also a parallel but vigorous
line of research concerning the possibility of defining the analogues of wave and scattering operators
for various classes of nonlinear evolution equations. 

\subsubsection{Lax--Phillips theory and radiation fields}
\label{sec:lax-phillips-theory}
In this final section we present yet another approach to scattering theory. This is the more abstract
approach developed by Lax and Phillips \cite{LP}, which has played an influential paradigmatic role. 
Directly following this we describe the theory pioneered by Friedlander \cite{Fr} on what he called the
radiation fields associated to solutions of a linear wave equation. These describe certain asymptotic
information about waves, and beyond their purely analytic appeal, they also provide a beautiful
realization of the Lax--Phillips theory. These radiation fields have received quite a lot of attention
in recent years, and the theory has been extended to various nonlinear settings as well. There is a 
forthcoming and much more detailed survey specifically about radiation fields \cite{MW} to which we 
direct the reader. 

Throughout this section we fix a Hilbert space $\mathcal{H}$ and a unitary semigroup $U(t)$ which
acts on it.  The specific application we have in mind is that $\mathcal{H}$ is the space $H_E$ of
 finite energy initial data for the wave equation on $\reals^{n}$ with $n$ odd and $U(t)$ is the wave
evolution operator.  More precisely, let $\mathcal{H}_{0}$ be the completion of the space
$\calC^{\infty}_{c}(\reals^{n})\times \calC^{\infty}_{c}(\reals^{n})$ with respect to the norm
\begin{equation*}
\norm[\mathcal{H}_{0}]{(\phi, \psi)}^{2} = \int_{\reals^{n}}\left(\left| \grad \phi (z)\right| ^{2} + \left|
\psi(z)\right|^{2}\right)\dz ;
\end{equation*}
then, for $(\phi, \psi)\in\mathcal{H}_{0}$, let $U_{0}(t)(\phi,\psi)$ be as defined in the previous section. 
The unitarity of $U_0$ corresponds to conservation of energy for solutions of this wave equation.

Return now to the general formulation.
\begin{defn}
\label{thm:1}
 A closed subspace $\mathcal{D} \subset \mathcal{H}$ is called \emph{outgoing}, respectively \emph{incoming}, if
\begin{enumerate}[(i)]
\item $U(t)\mathcal{D} \subset \mathcal{D}$ for $t > 0$, respectively $t<0$,  
\item $\bigcap_{t \in \reals}U(t)\mathcal{D} = \{ 0\}$, and
\item $\overline{\bigcup_{t\in\reals}U(t)\mathcal{D}} = \mathcal{H}$.
\end{enumerate}
\end{defn}

In the example above, the space $\mathcal{D}_{+}$ consists of the pairs $(\phi, \psi) \in \mathcal{H}_{0}$ 
for which the solution $u(t,z)$ vanishes for $|z| \leq t$ when $t \geq 0$.  Continuous
dependence of solutions of the wave equation on initial data shows that $\calD_+$ is a 
closed subspace.  The first and second properties follow from the observation that if 
$(\phi, \psi) \in \mathcal{H}_{0}$, then by finite propagation speed, the solution of the wave 
equation with initial data $U(s)(\phi, \psi)$ vanishes for $|z| \leq t+s$.  

The fourth property is more subtle.  For the unperturbed wave equation in odd dimensions, it is a consequence
of Huygens' principle; in even dimensions, one may prove it using local energy decay, but it can also be proved
fairly explicitly via the Radon transform. We say more about this later. 

The fundamental result of Lax--Phillips theory is the existence of a translation representation:
\begin{theorem}[Lax--Phillips \cite{LP}]
\label{thm:existence-of-translation-reps}
Let $U(t)$ be a group of unitary operators on $\mathcal{H}$, and $\mathcal{D}$ an outgoing subspace 
with respect to $U(t)$. Then there exists a Hilbert space $\mathcal{K}$ and an isometric isomorphism
\begin{equation*}
\Phi: \mathcal{H} \to L^{2}\left( (-\infty, \infty); \mathcal{K}\right)
\end{equation*}
such that $\Phi(\mathcal{D}) = L^{2}\left( (0,\infty); \mathcal{K}\right)$ and $\Phi \circ U(t) = T_{t}\circ \Phi$,
where $(T(t)f)(s) = f(s-t)$ is the standard translation action of $\RR$ on $L^2(\RR; \mathcal{K})$. 
The isomorphism $\Phi$ is unique up to an isomorphism of $\mathcal{K}$.  
\end{theorem}
The isomorphism given here is called an   \emph{outgoing translation representation} of $U(t)$.  
There is an essentially identical result giving an isomorphism $\Phi'$ which maps an incoming subspace
$\mathcal{D}_{-}$ to $L^{2}\left( (-\infty, 0); \mathcal{K}\right)$ and intertwines $U(t)$ with $T(t)$. 
This is called an \emph{incoming translation representation}. The auxilliary Hilbert
space $\mathcal{K}$ may be taken to be the same as for the outgoing translation representation,
but of course the map $\Phi'$ is different than $\Phi$. 

Returning again to the unperturbed wave equation in $\RR^n$, $n$ odd, there is an explicit way
to obtain the translation representations using the Radon transform.
\begin{defn}
\label{defn:radon-transform}
For any $f \in \calC^{\infty}_{c}(\reals^{n})$, define the \emph{Radon transform} 
\begin{equation*}
(Rf)(s, \theta) = \int_{\langle z, \theta\rangle = s}f(z)\dsigma (z),
\end{equation*}
where $\dsigma(z)$ is surface measure on the hyperplane $\langle z, \theta \rangle = s$.
Clearly $Rf \in \calC^\infty_c(\reals \times \sphere^{n-1})$. 
\end{defn}

A key property of the Radon transform for our purposes is that it is invertible and in fact
the inversion formula is quite explicit: 
\begin{equation*}
f(z) = \frac{1}{2\left( 2\pi\right)^{n-1}}\int_{\sphere^{n-1}}\left(\left| \PD[s]\right|^{n-1} Rf\right)(z\cdot \theta,
\theta)\dtheta, 
\end{equation*}
where $\left| \PD[s]\right|$ is defined by conjugating multiplication by $|\sigma|$ with respect to the
Fourier transform. A remarkable fact, which can be proved by direct computation, is that $R$ intertwines 
the Laplacians on $\RR^n$ and $\RR$, 
\begin{equation*}
R \lap f = \pd[s]^{2}Rf .
\end{equation*}

We now define the \emph{Lax--Phillips transform}: for $n$ odd, and $(\phi, \psi) \in \calC^{\infty}_{c}(\reals^{n}) \times 
\calC^{\infty}_{c}(\reals^{n})$, let
\begin{equation*}
\LP (\phi, \psi) (s,\theta) = \frac{1}{(2\pi) ^{(n-1)/2}} \left(
\PD[s]^{(n+1)/2}\left(R\phi\right) (s,\theta) - \PD[s]^{(n-1)/2}\left(R\psi\right)(s,\theta) \right) .
\end{equation*}
\begin{theorem}
\label{thm:LP-is-trans-rep}
For $n$ odd, the Lax--Phillips transform $\LP$ extends to a unitary isomorphism
\begin{equation*}
\LP : \mathcal{H}_{0} \to L^{2}\left( \reals ; L^{2}(\sphere^{n-1})\right),
\end{equation*}
and is a translation representation, 
\begin{equation*}
\left(\LP U_{0}(t) (\phi, \psi)\right) (s,\theta)= \left(T_{t}\LP (\phi, \psi) \right) (s,\theta) = \left( \LP (\phi,
\psi)\right)(s-t, \theta).
\end{equation*}
\end{theorem}

One consequence of Theorem~\ref{thm:LP-is-trans-rep} is that $\mathcal{H}_{0}$ splits as an orthogonal direct sum of
the incoming and outgoing subspaces:
\begin{equation}
\label{eq:split}
\mathcal{H}_{0} = \mathcal{D}_{+}\oplus \mathcal{D}_{-}.
\end{equation}
In particular, in this special case, the outgoing and incoming isomorphisms $\Phi$ and $\Phi'$ are equal. 

Now consider the wave equation with potential.  As before, assume that $n$ is odd and $V \in \calC^{\infty}_{c}(\reals^{n})$
is real-valued.  Choose $R$ so that $\supp V \subseteq B(0,R)$.  Let $U(t)$ be the group associated to the Cauchy problem
\begin{equation}
\label{eq:potential-eqn}
\Box u + Vu = 0, \qquad (u, \pd[t]u)|_{t=0} = (\phi, \psi),
\end{equation}
i.e.\ $U(t)(\phi, \psi) = (u(t), \pd[t]u(t))$. Since $V$ does not depend on $t$, there is a conserved 
energy, 
\begin{equation}
\label{eq:energy-norm}
\norm[E]{(u(t), \pd[t]u(t))}^{2} = \int_{\reals^{n}}\left( \left|\pd[t]u(t,z)\right|^{2} + \left| \grad u(t,z) \right|^{2} + V(z)
\left| u(t,z)\right|^{2}\right) \dz.
\end{equation}
The Hilbert space $\mathcal{H}$ is the set of pairs $(\phi, \psi)$ for which this energy is finite. It is not hard
to see, using the Sobolev inequality, that $\mathcal{H}$ and $\mathcal{H}_{0}$ consist of the same pairs of
elements, although the norm is different. The energy extends to the bilinear pairing on $\mathcal{H}$:
\begin{equation}
\label{eq:pairing}
  \left\langle \left(
    \begin{array}{c}
      \phi_{1} \\ \psi _{1}
    \end{array}
  \right) , \left(
    \begin{array}{c}
      \phi_{2}\\ \psi_{2}
    \end{array}
  \right) \right\rangle = \int_{\reals^{n}}\left( \grad\phi_{1} \cdot
  \grad \overline{\phi_{2}} + V(z)\phi_{1}\overline{\phi_{2}} +
  \psi_{1}\overline{\psi_{2}}\right) \dz .
\end{equation}

Consider now the operator
\begin{equation*}
A = \left(
    \begin{array}{cc}
      0 & 1 \\ \lap - V & 0 
    \end{array}\right);
\end{equation*}
this is anti-symmetric with respect to the pairing~\eqref{eq:pairing}.  The wave group $U(t)$ can be regarded
instead as the solution operator for the system
\begin{gather*}
\pd[t]
\begin{pmatrix}  u_{0} \\ u_{1} \end{pmatrix} = A
\begin{pmatrix}  u_{0} \\ u_{1} \end{pmatrix}, \qquad
\left.\begin{pmatrix}  u_{0} \\ u_{1} \end{pmatrix} \right|_{t=0} = 
\begin{pmatrix}  \phi \\ \psi \end{pmatrix}.
\end{gather*}

We now make a simplifying assumption that $L = -\lap + V$ has no $L^2$ eigenvalues, or equivalently, 
that $A$ has no such eigenvalues. Without this assumption, the results below require a projection
off the finite dimensional space $\calH_{\mathrm{pp}}$. We refer to \cite{Lax:1966} for more details
about how to proceed without this assumption. The advantage of this assumption is that now 
the energy~\eqref{eq:energy-norm} is positive definite. 

For this perturbed problem, we define the incoming and outgoing subspaces $\mathcal{D}_{\pm, R} \subset
\mathcal{H}$ to consist of those elements $(\phi, \psi)$ so that $U_{0}(t)(\phi, \psi)$ vanishes in $|z| \leq t
+ R$ for $t\geq 0$, respectively $|z| \leq - t + R$ for $t \leq 0$. Thus, in terms of the free incoming
and outgoing subspaces, $\mathcal{D}_{\pm, R} = U_{0}(\pm R)\mathcal{D}_{\pm}$.  The verification that
these satisfy all the correct properties relies on the following 
\begin{lemma}
\label{lem:same}
If $\mathbf{f}= (\phi, \psi)\in \mathcal{D}_{+,R}$, then $U_{0}(t)\mathbf{f} = U(t)\mathbf{f}$ for $t > 0$;
the analogous statement holds for $\mathcal{D}_{-,R}$ when $t < 0$.
\end{lemma}

We now use this to show that $\mathcal{D}_{+,R}$ is an outgoing subspace for $U(t)$ on $\mathcal{H}$.  
Indeed, by this lemma, the first two properties follow from the corresponding properties of $\calD_+$.
For the third property, suppose we know that for any compact subset $K\subset \reals^{n}$ and any 
solution $u$ of \eqref{eq:potential-eqn}, we have
\begin{multline*}
\lim_{t\to \infty}\norm[E,K]{u(t)}^{2} := \\
\lim_{t\to \infty} \int_{K}\left( \left| \pd[t]u(t,z)\right|^{2} + 
\left|\grad u(t,z)\right|^{2} + V(z) |u(t,z)|^{2} \right) \dz = 0 .
\end{multline*}
This is called local energy decay, and is known to be true in many circumstances.  
Now consider the initial data $\mathbf{f} = (\phi, \psi) \in
\mathcal{H}$ with $\mathbf{f}\perp \bigcup U(t)\mathcal{D}_{+,R}$ with respect to the 
pairing~\eqref{eq:pairing}. Thus $U(t) \mathbf{f} \perp \mathcal{D}_{+,R}$ for any $t$, and
in particular, $U(t)\mathbf{f} \perp \mathcal{D}_{+,R}$ with respect to the standard pairing on
$\dot{H}^{1}\times L^{2}$.  This shows that $U_{0}(-R)U(t)\mathbf{f}\perp \mathcal{D}_{+}$ with 
respect to the standard pairing, and hence $U_{0}(-R)U(t)\mathbf{f}\in \mathcal{D}_{-}$ and
$U_{0}(-2R)U(t)\mathbf{f}\in \mathcal{D}_{-,R}$.  

Consider now $v(s,z) = U(s)U_{0}(-2R)U(t)\mathbf{f}$. By Lemma~\ref{lem:same}, $v(s,z)$ agrees with
$U_{0}(s)U_{0}(-2R)U(t)\mathbf{f}$ for $s < 0$ and thus vanishes for $|z|\leq -s+R$ for $s<0$.  

Now we bring in the local energy decay. This implies that for any $\epsilon > 0$, if $t$ 
is sufficiently large then $\norm[E,B(5R)]{U(t)\mathbf{f}} < \epsilon$. For such $t$, 
finite propagation speed  implies both
\begin{equation*}
  \norm[E,B(3R)]{U_{0}(-2R)U(t)\mathbf{f}} < \epsilon, \quad
  \text{and}\quad \norm[E,B(3R)]{U(-2R)U(t)\mathbf{f}} < \epsilon .
\end{equation*}
Because the two equations and the initial data agree outside $B(R)$, using finite propagation speed
again, we get that $U_{0}(-2R)U(t)\mathbf{f} =U(-2R)U(t)\mathbf{f}$ for $|z| > 3R$ and hence 
\begin{equation*}
  \norm[E]{U_{0}(-2R)U(t)\mathbf{f} - U(-2R)U(t)\mathbf{f}}< 2\epsilon .
\end{equation*}
Because $U(t)$ is unitary with respect to \eqref{eq:pairing}, applying $U(2R-t)$ to the difference shows that
\begin{equation*}
  \norm[E]{U(2R-t)U_{0}(-2R)U(t)\mathbf{f} - \mathbf{f}} < 2\epsilon .
\end{equation*}
Finally, since $t$ is large, $2R-t < 0$ and so
$U(2R-t)U_{0}(-2R)U(t)\mathbf{f} = U_{0}(2R-t)U_{0}(-2R)U(t)\mathbf{f}$ by Lemma~\ref{lem:same}.  This
shows that in fact
\begin{equation*}
  \norm[E]{U_{0}(-t)U(t)\mathbf{f} - \mathbf{f}}<2\epsilon.
\end{equation*}
Because $U_{0}(-R)U(t)\mathbf{f}\in \mathcal{D}_{-}$, the first term here is an element of $\mathcal{D}_{-, t-R}$ 
and thus vanishes for $|z| \leq t-R$.  Taking $t$ even larger gives 
\begin{equation*}
\norm[E]{\mathbf{f}} < 2\epsilon,
\end{equation*}
and therefore $\mathbf{f} = 0$.  This establishes the third property.

\medskip

Theorem~\ref{thm:existence-of-translation-reps} asserts the existence of incoming and outgoing translation 
representations for the incoming and outgoing subspaces $\mathcal{D}_{-,R}$ and $\mathcal{D}_{+,R}$.  
We shall give a a concrete realization of these using the so-called radiation fields. 

Our next goal is to show that a particular quantitative rate of local energy decay implies that the local
energy actually decays exponentially.  
\begin{theorem}
\label{thm:exp-decay}
Suppose that for each compact subset $K \subset \reals^{n}$ there is a function $c_{K}(t)$ which tends to $0$ as 
$t\to \infty$, such that if the Cauchy data $u(0)$ have support in $K$, then 
\begin{equation}
\label{eq:quant-decay}
\norm[E,K]{u(t)}^{2} \leq c_{K}(t) \norm[E]{u(0)}^{2}.
\end{equation}
Then there are positive constants $C$ and $\alpha$ depending on $K$ such that if $u(0)$ is supported in $K$, then
\begin{equation}
\label{eq:exp-decay}
\norm[E,K]{u(t)} \leq Ce^{-\alpha t} \norm[E]{u(0)}
\end{equation}
for all $t > 0$.  
\end{theorem}

The proof uses the compactness properties of the Lax--Phillips semigroup $Z(t)$, which we introduce now.  If $
P_{\pm,R}$ are the orthogonal projections onto the orthocomplements of $\mathcal{D}_{\pm,R}$, then $Z(t)$ is given 
for $t\geq 0$ by 
\begin{equation*}
  Z(t) = P_{+,R}U(t)P_{-,R}.
\end{equation*}
The local energy decay hypothesis in the theorem statement implies that, for $t$ large enough, $Z(t)$ has norm 
bounded by $1/2$, and repeated application of $Z(t)$ leads to the exponential decay.

We are now in a position to introduce the radiation field of a solution $u$ to the perturbed wave equation.  
The idea is to identify initial data for $u$  with a normalized limit of the solution along outgoing (or incoming) 
light rays. As before, we start with the definition of these radiation fields for the unperturbed operator. 

Suppose that $u$ solves $\Box_{0}u =0$ with intial data $(\phi,
\psi)$. Introduce coordinates $s = t-|z|$ and $x = |z|^{-1}$; these
parametrize the family of outgoing light rays and the position along
them.  Now define the auxiliary function
\[
v : \reals_{s} \times (0, \infty)_{x} \times \sphere^{n-1}_{\theta}  \to \reals, \qquad 
v(s, x, \theta) = x^{-\frac{n-1}{2}}u \left(s + \frac{1}{x}, \frac{1}{x}\theta \right)
\]
Here $\frac{1}{x}\theta$ is simply $z$ in polar coordinates. Since $x^{2}g_{M}$ is nondegenerate at $x=0$, 
the function $v$ extends smoothly across $x=0$. We then define the forward radiation field operator $\mathcal{R}_{+}$:
\begin{equation*}
  \mathcal{R}_{+}  (\phi, \psi)(s, \theta) = \pd[s]v(s,0,\theta)
\end{equation*}
The derivative of $v$ is included here to make $\mathcal R_+$ an isometric isomorphism:
\begin{equation*}
\mathcal{R}_{+}: \mathcal{H}_{0} \to L^{2}(\reals \times \sphere^{n-1})
\end{equation*}
Furthermore, the Minkowski metric is static, so $\mathcal R_+$  intertwines wave evolution and translation in $s$:
\begin{equation*}
\mathcal{R}_{+}U_{0}(T)(\phi, \psi) (s,\theta) = \mathcal{R}_{+}(\phi, \psi) (s-T, \theta)
\end{equation*}

Now observe that if $\mathbf{f}\in\mathcal{D}_{+}$, then $\mathcal{R}_{+}\mathbf{f}$ vanishes when $s \geq 0$.  
This follows from the unitarity of the radiation field operator, and the inverse image of those functions in 
$L^{2}(\reals\times \sphere^{n-1})$ supported in the nonpositive half-cylinder form an outgoing subspace
$\widetilde{\mathcal{D}}_{+}$.
\begin{equation*}
\widetilde{\mathcal{D}}_{+} = \left\{ \mathcal{R}^{-1}_{+}f : f(s,\theta) = 0    \text{ for }s > 0\right\}
\end{equation*}
Indeed, this is a closed subspace; the first and second  properties follow directly from the fact that $\mathcal{R}_{+}$ is 
a translation representation, while the third property follows from the surjectivity of $\mathcal{R}_{+}$.  
One may also define $\widetilde{\mathcal{D}}_{-}$ via the backward radiation field $\mathcal{R}_{-}$; this encodes 
information from solutions in the limit as $t \searrow -\infty$.   For the free wave equation, $\mathcal R_+$ 
has an explicit expression in terms of the Radon transform, and this can be used to show that $\mathcal{H}_{0} =
\widetilde{\mathcal{D}}_{+}\oplus\widetilde{\mathcal{D}}_{-}$.  

For the perturbed equation the forward and backward radiation fields, $\mathcal{L}_{\pm}$, are defined
in the same way. We can also define the scattering operator $\mathcal{A}$ using the radiation fields by
\begin{equation*}
\mathcal{A} = \mathcal{L}_{+}\mathcal{L}_{-}^{-1}
\end{equation*}
Thus $\mathcal A$ maps data at past null infinity into data at future null infinity.  The relationship
to the scattering operator $\mathcal{S}$ introduced in Section~\ref{sec:moll-wave-oper} is that
\begin{equation*}
\mathcal{S} = \mathcal{R}_{+}^{-1}\mathcal{L}_{+}\mathcal{L}^{-1}_{-}\mathcal{R}_{-}
= \mathcal{R}_{+}^{-1}\mathcal{A}\mathcal{R}_{-}
\end{equation*}
The conjugation of $\mathcal A$ by the Fourier transform in $s$ corresponds to the \emph{scattering matrix} 
employed by Melrose in \cite{Melrose:1994}. 

The radiation field exists and is a unitary operator in a variety of geometric settings.  On asymptotically Euclidean 
spaces, this is due to Friedlander~\cite{Fr,Friedlander:2001} and S{\'a} 
Barreto~\cite{Sa-Barreto:2003,Sa-Barreto:2008}; on asymptotically real and complex hyperbolic manifolds it was proved
by S{\'a} Barreto~\cite{Sa-Barreto:2005}, and Guillarmou and S{\'a} Barreto~\cite{Guillarmou:2008}, respectively.  
In the asymptotically Euclidean and real hyperbolic cases, S{\'a} Barreto and Wunsch~\cite{SBW}
prove that it is a Fourier integral operator with canonical relation given by the sojourn relation, a close relative of 
the Busemann function in each of these geometric settings. 
The radiation field has also been defined in certain nonlinear and non-static situations.  In particular, the
first author and S{\'a} Barreto show \cite{BSB} that it exists and is norm-preserving for certain semilinear wave 
equations in $\reals^{3}$, while Wang~\cite{Wang:2010,Wang:2011} defined the radiation field for the 
Einstein equations on perturbations of Minkowski space when the spatial dimension is at least $4$.  
Forthcoming work of the first author, Vasy, and Wunsch~\cite{BVW} analyzes the $s \to \infty$ asymptotics 
of the radiation field on (typically non-static) perturbations of Minkowski space.


\end{document}